\documentclass[preprint,12pt,1p]{elsarticle}
\makeatletter
\def\ps@pprintTitle{%
  \let\@oddhead\@empty
  \let\@evenhead\@empty
  \let\@oddfoot\@empty
  \let\@evenfoot\@oddfoot
}
\makeatother
\usepackage{graphics}
\usepackage{float}
\usepackage{placeins}
\usepackage{adjustbox}
\usepackage{epsfig}
\usepackage{mathptmx}
\usepackage{amsfonts}
\usepackage{amsmath}
\usepackage{amssymb}
\usepackage{hyperref}
\usepackage{fullpage}
\usepackage{amsthm}
\usepackage{geometry}
\theoremstyle{plain}

\theoremstyle{definition}

\theoremstyle{example}

\theoremstyle{remark}

\newlength{\defbaselineskip}
\journal{Computer Aided Geometric Design}
\numberwithin{equation}{section}
\begin{document}
\begin{frontmatter}
\title{ A Generalized Closed Form of Ramanujan-Type Fourier Cosine
Transform via Meijer’s $G$-Function}
\author[1]{*S. A. Dar}
\ead{showkatjmi34@gmail.com}
\author[2]{R. P. Paris}

\address[1]{Department of Mathematics and Statistics, Govt Degree College Sopore, \\Jammu and Kashmir Higher Education Department,  India\\E-Mail: showkatjmi34@gmail.com}
\address[2]{Division of Computing and Mathematics, Abertay University, Dundee DD1 1HG, UK\\ E-Mail: r.paris@abertay.ac.uk}
\cortext[cor1]{Corresponding author}
\begin{abstract}
In this paper, we obtain analytical evaluations of the Ramanujan integral
\[\textbf{R}_{C}(m,n)= \int_{0}^{\infty}\frac{x^m\,\cos(\pi nx)}{\exp{(2\pi\sqrt{x})-1}}dx\]
subject to suitable convergence conditions in terms of an infinite series of Meijer $G$-functions of one variable, by using Mellin-Barnes-type contour integral representations of the cosine function. 
We also consider some generalizations of the integral $\textbf{R}_{C}(m,n)$ given as the integrals $I_{C}^{*}(\upsilon,b,c,\lambda,y)$ ,$\Xi_{C}(\upsilon,b,c,\lambda,y)$, $\nabla_{C}(\upsilon,b,c,\lambda,y)$ and $I_{C}(\upsilon,b,\lambda,y)$.  These integrals are also expressed in terms of infinite series of Meijer $G$-functions. Moreover, as an application of a Ramanujan's integral $\textbf{R}_C(m,n)$, the closed-form evaluations of nine infinite series of Meijer $G$-functions are obtained.\\
\\
\textit{2010 AMS Classification: 33C60, 33C65, 33C70, 33B15, 33C05, 33C45, 33C10}
 \end{abstract}
\end{frontmatter}
\section{Introduction and Preliminaries}
\label{Sec:1}

The Fourier cosine transform of $g(x)$ over the interval $[0,\infty)$ is defined by\\
\begin{equation}\label{RCG30}
F_{C}\{g(x);y\}=\int_{0}^{\infty}g(x)\cos(xy)dx=G_{C}(y),~~~~~(y>0).
\end{equation}
 In the literature  \cite{C4,D,E2,o,s3}, the analytical evaluations of the Fourier cosine transform of $x^{\upsilon-1}/(\exp (bx)\pm 1)$
are available in terms of Riemann's zeta function, the Psi function (or Digamma function) and the hyperbolic and Beta functions.
In \cite[p. 85, Eq.~(49), last line]{R1}, Ramanujan considered the evaluation of the integral
\begin{equation}\label{RCG1}
\textbf{R}_{C}(m,n)= \int_{0}^{\infty}\frac{x^m\,\cos(\pi nx)}{\exp{(2\pi\sqrt{x})-1}}dx,
\end{equation}
and showed that \cite[p.86, Eq.~(50)]{R1}:
\begin{equation}\label{RCG2}
\textbf{R}_{C}(1,1/2)=\int_{0}^{\infty}\frac{x\cos(\frac{\pi x}{2})}{\exp{(2\pi\sqrt{x})-1}}dx=\frac{13-4\pi}{8\pi^{2}},~~~~~~~~~~~~~~~~
\end{equation}
\begin{equation}\label{RCG3}
\textbf{R}_{C}(1,2)=\int_{0}^{\infty}\frac{x\cos(2\pi x)}{\exp{(2\pi\sqrt{x})-1}}dx=\frac{1}{64}\left(\frac{1}{2}-\frac{3}{\pi}+\frac{5}{\pi^{2}}\right),
\end{equation}
\begin{equation}\label{RCG4}
\textbf{R}_{C}(2,2)=\int_{0}^{\infty}\frac{x^{2}\cos(2\pi x)}{\exp{(2\pi\sqrt{x})-1}}dx=\frac{1}{256}\left(1- \frac{5}{\pi}+\frac{5}{\pi^{2}}\right).
\end{equation}
When $m=0$ and for particular values of $n$, some values of $\textbf{R}_{C}(0,n)$ \cite[p.85 (eq. 48)]{R1} are given below:
\begin{equation}\label{RCG9}
\textbf{R}_{C}(0,0)=\int_{0}^{\infty}\frac{dx}{\exp{(2\pi\sqrt{x})-1}}=\frac{1}{12},
\end{equation}
\begin{equation}\label{RCG10}
\textbf{R}_{C}(0,1)=\int_{0}^{\infty}\frac{\cos(\pi x)}{\exp{(2\pi\sqrt{x})-1}}dx=\frac{2-\sqrt{2}}{8},
\end{equation}
\begin{equation}\label{RCG11}
\textbf{R}_{C}(0,2)=\int_{0}^{\infty}\frac{\cos(2\pi x)}{\exp{(2\pi\sqrt{x})-1}}dx=\frac{1}{16},
\end{equation}
\begin{equation}\label{RCG12}
\textbf{R}_{C}(0,4)=\int_{0}^{\infty}\frac{\cos(4\pi x)}{\exp{(2\pi\sqrt{x})-1}}dx=\frac{3-\sqrt{2}}{32},
\end{equation}
\begin{equation}\label{RCG13}
\textbf{R}_{C}(0,6)=\int_{0}^{\infty}\frac{\cos(6\pi x)}{\exp{(2\pi\sqrt{x})-1}}dx=\frac{13-4\sqrt{3}}{144},
\end{equation}
\begin{equation}\label{RCG14}
\textbf{R}_{C}(0,1/2)=\int_{0}^{\infty}\frac{\cos\left(\frac{\pi x}{2}\right)}{\exp{(2\pi\sqrt{x})-1}}dx=\frac{1}{4\pi},
\end{equation}
\begin{equation}\label{RCG15}
\textbf{R}_{C}(0,2/5)=\int_{0}^{\infty}\frac{\cos\left(\frac{2\pi x}{5}\right)}{\exp{(2\pi\sqrt{x})-1}}dx=\frac{8-3\sqrt{5}}{16}.
\end{equation}

The connection between $\textbf{R}_C(0,n)$ and the corresponding Fourier sine integral is given by the following theorem \cite[p.76-77, Eqs.~10, $10^{'}$)]{R1}.
If
\begin{equation}\label{RCG6}
\Upsilon(n)=\frac{1}{2\pi n}+\int_{0}^{\infty}\frac{\sin(\pi nx)}{\exp{(2\pi\sqrt{x})-1}}dx,
\end{equation}
then
\begin{equation}\label{RCG7}
\textbf{R}_{C}(0,n)=\frac{1}{n}\sqrt{{\frac{2}{n}}}\,\Upsilon\left(\frac{1}{n}\right)-\Upsilon(n),
\end{equation}
and
\begin{equation}\label{RCG8}
\Upsilon(n)=\frac{1}{n}\sqrt{\frac{2}{n}}\,\textbf{R}_C\left(0,\frac{1}{n}\right)+\textbf{R}_C(0,n),
\end{equation}
where $n$ is positive rational number.

A natural generalization of the Gauss hypergeometric function ${}_{2}F_{1}(z)$ is the generalized hypergeometric function ${}_{p}F_{q}(z)$  
 with $p$ numerator parameters $\alpha_{1},  ... , \alpha_{p}$ and $q$ denominator parameters $\beta_{1}, ... , \beta_{q}$ defined by \cite{R}\\
\begin{equation}\label{FAR4}
{}_{p}F_{q}\left(\ \begin{array}{lll}\alpha_{1}, ... ,\alpha_{p}~;~\\\beta_{1}, ... , \beta_{q}~;~\end{array} z\right)
=\sum_{n=0}^{\infty}\frac{(\alpha_{1})_{n}  ...  (\alpha_{p})_{n}}{(\beta_{1})_{n} ... (\beta_{q})_{n}}\frac{z^{n}}{n!}~,
\end{equation}
where $\alpha_{j}\in\mathbb{C}~(j=1, ... ,p)$ , $\beta_{j}\in\mathbb{C}\setminus \mathbb{Z}_{0}^{-}~(j=1, ... ,q)$ and
$\ p,~q\in\mathbb{N}_{0}$.
Then the hypergeometric ${}_{p}F_{q}(z)$ function in (\ref{FAR4}) converges absolutely for $|z|<\infty$ when $p\leq q$ and for $|z|<1$ when $p=q+1$. Furthermore, if we set,
\begin{equation}\label{FAR6}
\omega:=\left(\ \sum_{j=1}^{q}\beta_{j}-\sum_{j=1}^{p}\alpha_{j}\right),
\end{equation}
it is known that when $p=q+1$ the function ${}_{p}F_{q}(z)$  is absolutely convergent for $|z|=1$ if $\Re(\omega)>0$, conditionally convergent for $|z|=1$ ($z\neq 1$) if $-1<\Re (\omega)<0$ and divergent for $|z|=1$ if $\Re (\omega)\leq -1$.
The important special cases of ${}_{p}F_{q}(z)$
include (for example) the binomial series ${}_{1}F_{0}(z)$  given by
\begin{equation}\label{RCG20}
(1-z)^{-a}={}_{1}F_{0}\left(\begin{array}{lll}a~;\\  \overline{~~~};\end{array} z \right)
=\sum_{n=0}^{\infty}\frac{(a)_{n}}{n!}z^{n},
\end{equation}
where $|z|<1$, $a\in\mathbb{C}$.

For every positive integer $m$ \cite[p.~22,(26)]{s1}, we have the Gauss-Legendre multiplication theorem for the Gamma function
\begin{equation}\label{RCG50}
  \Gamma(mz)=(2\pi)^{\frac{(1-m)}{2}}m^{mz-\frac{1}{2}}\prod_{j=1}^{m}\Gamma\left(z+\frac{j-1}{m}\right),~~ mz\in\mathbb{C}\setminus \mathbb{Z}_{0}^{-}.
\end{equation}
The Fox-Wright function ${}_{p}\Psi_{q}(z)$ of one variable \cite{K5,K6} is given by
\begin{eqnarray}\label{RCG21}
{}_{p}\Psi_{q}\left[\ \begin{array}{lll}(\alpha_{1}, A_{1}),...,(\alpha_{p}, A_{p});\\(\beta_{1}, B_{1}),...,(\beta_{q}, B_{q});~\end{array} z \right]
=\sum_{k=0}^{\infty}\frac{\Gamma(\alpha_{1}+kA_{1})...\Gamma(\alpha_{p}+kA_{p})}{\Gamma(\beta_{1}+kB_{1})...\Gamma(\beta_{q}+kB_{q})}\frac{z^{k}}{k!},
\end{eqnarray}
\begin{eqnarray}\label{RCG25}
~~~~~~~~~~~~~~~~~~~~~~~~~~~~~~~~~~~~~~~~~~~~~~~~~~~~~~~~~
=\frac{1}{2\pi i}\int_{L}\frac{\displaystyle \Gamma(-\zeta)\prod_{i=1}^{p}\Gamma(\alpha_{i}+A_{i}\zeta)}{\displaystyle\prod_{j=1}^{q}\Gamma(\beta_{j}+B_{j}\zeta)}(-z)^{\zeta}d\zeta~,~~~~~~~~~~~~~~
\end{eqnarray}
where $z\in\mathbb{C};~$ the parameters $\alpha_{i},\beta_{j}\in\mathbb{C}$ and coefficients $A_{i},B_{j}\in\mathbb{R}=(-\infty,+\infty)$ in the case of series (\ref{RCG21})  (or $A_{i},B_{j}\in\mathbb{R}_{+}=(0,+\infty)$ in the case of the contour integral (\ref{RCG25})), $A_{i}\neq0~(i=1,2,...,p),B_{j}\neq0~(j=1,2,...,q)$. In the series (\ref{RCG21}), the parameters $\alpha_{i},\beta_{j}$ and coefficients $A_{i},B_{j}$ are adjusted in such a way that the products of Gamma functions in the numerator and denominator are well defined.

Suppose:
\begin{equation}\label{RCG26}
\Delta^{*}=\left(\ \sum_{j=1}^{q}B_{j}-\sum_{i=1}^{p}A_{i}\right),
\end{equation}
\begin{equation}\label{RCG27}
\delta^{*}=\left(\prod_{i=1}^{p}|A_{i}|^{-A_{i}}\right)\left(\prod_{j=1}^{q}|B_{j}|^{B_{j}}\right),
\end{equation}
\begin{equation}\label{RCG28}
 \mu^{*}=\sum_{j=1}^{q}\beta_{j}-\sum_{i=1}^{p}\alpha_{i}+\left(\frac{p-q}{2}\right),
\end{equation}
and
\begin{equation}\label{RCG29}
\sigma^{*}=(1+A_{1}+...+A_{p})-(B_{1}+...+B_{q})=1-\Delta^{*}.
\end{equation}
Then we have the following convergence conditions for (\ref{RCG21}) and (\ref{RCG25}) when the contour $L$ starts from $\gamma-i\infty$ and ends at $\gamma+i\infty$, where $\gamma\in\mathbb{R}$:
(i) When $\sigma^{*}>0$, $|\arg(-z)|<\frac{\pi}{2}\sigma^{*}$,  $0<|z|<\infty~,z\neq0$;\\
(ii)  When $\sigma^{*}=0$,  $\arg(-z)=0$ , $0<|z|<\infty,~z\neq0$ such that $-\gamma \Delta^{*}+Re(\mu^{*})>\frac{1}{2}+\gamma$;\\
(iii) When $\gamma=0$, $\sigma^{*}=0, \arg(-z)=0$, $0<|z|<\infty,~z\neq0$, such that $Re(\mu^{*})>\frac{1}{2}$.

The Wright function is defined by [p. 438, Eq.(1.2)]\cite{K6}
\begin{equation}
\Phi(\alpha,\beta,z)={}_{0}\Psi_{1}\left[\ \left.\begin{array}{lll}\overline{~~~~~~~~~}\\(\beta, \alpha)\end{array}\right| z\right]=\sum_{k=0}^{\infty}{\frac{1} {\Gamma(\beta+\alpha k)}\frac{z^{k}}{k!}},
\end{equation}
where $z,\beta\in\mathbb{C}$ ; $\alpha>0$.
The Wright generalized Bessel function [p. 438, Eq.(1.3)]\cite{K6}, is defined by
\begin{equation}
J_{\nu}^{\mu}(-z)=\Phi(\mu,\nu+1,-z)={}_{0}\Psi_{1}\left[\ \left.\begin{array}{lll}\overline{~~~~~~~~~}\\(\nu+1, \mu)\end{array}\right| -z\right]=
\sum_{k=0}^{\infty}\frac{1} {\Gamma(\nu+1+\mu k)}\frac{(-z)^{k}}{k!},
\end{equation}
where $z,\nu\in\mathbb{C}$ and $\mu>0$.
The two-parameter Mittag-Leffler function [p.450, Eq.(6.1)]\cite{K6}, is given by
\begin{equation}
E_{\alpha,\beta}(z)={}_{0}\Psi_{1}\left[\ \left.\begin{array}{lll}(1,1)  \\(\beta, \alpha)\end{array}\right| z\right]=\sum_{k=0}^{\infty}\frac{z^k} {\Gamma(\beta+\alpha k)},
\end{equation}
where $z,\beta\in\mathbb{C}$ and $\alpha>0$.

The Meijer $G$-function is defined by means of  the Mellin-Barnes contour integral \cite[Sec.(1.5), (1)]{s1}.
\begin{multline}\label{RCG32}
G_{{p} ,{q}}^{{m},{n}} \left( z\right):=
G_{{p} ,{q}}^{{m},{n}} \left( z~\bigg{|} \begin{array}{lll} \alpha_{p} \\ \beta{q} \end{array} \right):=
G_{{p} ,{q}}^{{m},{n}} \left( z~\bigg{|} \begin{array}{lll} \alpha_{1},...,\alpha_{n};\alpha_{n+1},...,\alpha_{p} \\ \beta_{1},...,\beta_{m};\beta_{m+1},...,\beta{q} \end{array} \right)\\
=\frac{1}{2\pi i}\int_{L}
\frac{\displaystyle\prod_{j=1}^{m}\Gamma(\beta_{j}-\zeta)\prod_{j=1}^{n}\Gamma(1-\alpha_{j}+\zeta)}{\displaystyle\prod_{j=m+1}^{q}\Gamma(1-\beta_{j}+\zeta)\prod_{j=n+1}^{p}\Gamma(\alpha_{j}-\zeta)}z^{\zeta}d\zeta~,~~~~~~~~~~~~~~
\end{multline}
where $L$ is an upward oriented contour which separates the poles of $\Gamma(\beta_{m}-\zeta)$ from those of $\Gamma(1-\alpha_{n}+\zeta)$ and which begins and ends at $-i\infty$ and $i\infty$. Throughout this work, we assume
 $z\neq0$, and $m,n,p,q$ are non-negative integers such that $ 1\leq m\leq q$~,~$0\leq n\leq p$ and $p\leq q$.
The integral (\ref{RCG32}) converges in the sector $|\arg(z)|<\pi \kappa$, where $\kappa=m+n-\frac{1}{2}(p+q)$ and it is supposed that $\kappa>0$.
%

The representations of cosine function in terms of  Meijer's $G$-Function \cite{E1, E2} are given by
\begin{equation}\label{RCG38}
\cos(z)= \sqrt{\pi}~~G_{{0} ,{2}}^{{1},{0}} \left(\ \frac{z^{2}}{4}\bigg{|} \begin{array}{lll}\overline{~~~~~~~~~~~~~}\\ 0~~;~~ \frac{1}{2}\end{array} \right),
\end{equation}
\begin{equation}\label{RCG43}
=\sqrt{\pi}\frac{1}{2\pi i}\int_{-i\infty}^{+i\infty}\frac{\Gamma(-\zeta)}{\Gamma(\frac{1}{2}+\zeta)}\left(\frac{z^{2}}{4}\right)^{\zeta}d\zeta,
\end{equation}
\begin{equation}\label{RCG39}
\cos(z)= \sqrt{\frac{\pi z}{2}}~~G_{{0} ,{2}}^{{1},{0}} \left(\ \frac{z^{2}}{4}\bigg{|} \begin{array}{lll}\overline{~~~~~~~~~~~~~}\\ \frac{-1}{4}~~;~~\frac{1}{4}\end{array} \right),
\end{equation}
\begin{equation}\label{RCG44}
~~~~~~=\sqrt{\frac{\pi z}{2}}\frac{1}{2\pi i}\int_{-i\infty}^{+i\infty}\frac{\Gamma\left(-\frac{1}{4}-\zeta\right)}{\Gamma\left(\frac{3}{4}+\zeta\right)}\left(\frac{z^{2}}{4}\right)^{\zeta}d\zeta,
\end{equation}
\begin{equation}\label{RCG40}
\cos(z)= \frac{2\sqrt{\pi}}{z}~~G_{{0} ,{2}}^{{1},{0}} \left(\ \frac{z^{2}}{4}\bigg{|} \begin{array}{lll}\overline{~~~~~~~~~~~~~}\\\frac{1}{2}~~;~~ 1\end{array} \right),
\end{equation}
\begin{equation}\label{RCG45}
=\frac{2\sqrt{\pi}}{z}\frac{1}{2\pi i}\int_{-i\infty}^{+i\infty}\frac{\Gamma(\frac{1}{2}-\zeta)}{\Gamma(\zeta)}\left(\frac{z^{2}}{4}\right)^{\zeta}d\zeta,~~~~
\end{equation}

\begin{equation}\label{RCG41}
\cos(z)= \sqrt{\pi}\left(\frac{z}{2}\right)~~G_{{0} ,{2}}^{{1},{0}} \left(\ \frac{z^{2}}{4}\bigg{|} \begin{array}{lll}\overline{~~~~~~~~~~~~~}\\-\frac{1}{2}~~;~~0\end{array} \right),
\end{equation}
\begin{equation}\label{RCG46}
~~~~~~~=\sqrt{\pi}\left(\frac{z}{2}\right)\frac{1}{2\pi i}\int_{-i\infty}^{+i\infty}\frac{\Gamma(-\frac{1}{2}-\zeta)}{\Gamma(1+\zeta)}\left(\frac{z^{2}}{4}\right)^{\zeta}d\zeta,
\end{equation}
\begin{equation}\label{RCG42}
\cos(z)= \sqrt{2\pi}~~G_{{0} ,{4}}^{{2},{0}} \left(\ \frac{z^{4}}{256}\bigg{|} \begin{array}{lll}\overline{~~~~~~~~~~~~~~~~~~~~~~~~}\\0~~,~~ \frac{1}{2}~~;~~\frac{1}{4}~~,~~\frac{3}{4}\end{array} \right),
\end{equation}
\begin{equation}\label{RCG47}
~~~~~~~~~~~~~~~~=\sqrt{2\pi}\frac{1}{2\pi i}\int_{-i\infty}^{+i\infty}\frac{\Gamma(-\zeta)\Gamma(\frac{1}{2}-\zeta)}{\Gamma(\frac{3}{4}+\zeta)\Gamma(\frac{1}{4}+\zeta)}\left(\frac{z^{4}}{256}\right)^{\zeta}d\zeta,
\end{equation}
where $z\neq0$.

In the available literature \cite{B5,B6,R1,R2} on Ramanujan's work, the analytical evaluation of Ramanujan's integral $\textbf{R}_{C}(m,n)$ is not given. Therefore, the main aim of this paper is to obtain representations of $\textbf{R}_{C}(m,n)$ in terms of Meijer's $G$-function of one variable. Also, our work is motivated by the papers given in \cite{M,Q12,Q2,Q3,R4}.
In this paper, we generalize Ramanujan's integral $\textbf{R}_{C}(m,n)$  in the following forms:
\[I_{C}^{*}(\upsilon,b,c,\lambda,y)=\displaystyle{\sum_{k=0}^{\infty}}\frac{\Theta(k)}{k!}
\int_{0}^{\infty}~x^{\upsilon-1}e^{-(\lambda b +c k)\sqrt{x}}\cos (xy)\,dx,\]
\[\Xi_{C}(\upsilon,b,c,\lambda,y)=\displaystyle{\int_{0}^{\infty}}x^{\upsilon-1}e^{-b\lambda\sqrt{x}}
{}_{r}\Psi_{s}\left[\ \begin{array}{lll}(\alpha_{1}, A_{1}),...,(\alpha_{r}, A_{r});\\(\beta_{1}, B_{1}),...,(\beta_{s}, B_{s});~\end{array} e^{-c\sqrt{x}}\right]\cos (xy)\,dx,\]
\[\nabla_{C}(\upsilon,b,c,\lambda,y)=\displaystyle{\int_{0}^{\infty}}x^{\upsilon-1}e^{-b\lambda\sqrt{x}}{}_{r}F_{s}\left(\ \begin{array}{lll}\alpha_{1},...,\alpha_{r};\beta_{1},...,\beta_{s};~\end{array} e^{-c\sqrt{x}}\right)\cos (xy)\,dx\]
and
\[I_{C}(\upsilon,b,\lambda,y)=
\displaystyle{\int_{0}^{\infty}}x^{\upsilon-1}\{\exp(b\sqrt{x})-1\}^{-\lambda}\cos(xy)\,dx,\]
where $\{\Theta(k)\}_{k=0}^{\infty}$ is a bounded sequence and obtain the analytical representations. Moreover, we also show how the main general theorems given below are applicable for obtaining new interesting results  by suitable adjustment of the parameters and variables.


\section{Generalized Ramanujan's integral in terms of Meijer $G$-functions}
\label{Sec:2}
We give some analytical expressions of the generalized Ramanujan integral $ I_{C}^{*}(\upsilon,b,c,\lambda,y)$ in the form of infinite series of the Meijer $G$-function\\
\textbf{Theorem 1.} {\it Let $\{\Theta(k)\}_{k=0}^{\infty}$ be a bounded sequence of arbitrary real or complex numbers
and suppose that the conditions $\Re(\upsilon)>0$, $\Re(c)>0$, $y>0$, $\Re(\lambda)>0$, $\Re(b)>0$ (or $\Re(\lambda)<0$,$\Re(b)<0$)}. Then we have
\begin{equation*}\label{RG1}
I_{C}^{*}(\upsilon,b,c,\lambda,y)=\sum_{k=0}^{\infty}\frac{\Theta(k)}{k!}
\int_{0}^{\infty}~x^{\upsilon-1}e^{-(\lambda b +c k)\sqrt{x}} \cos (xy)\,dx~~~~~~~~~~~~~~~~~~~~~~~~~~~~~~~~~~~~~~~~~
\end{equation*}
\begin{eqnarray}\label{RG2}
=\frac{2^{4\upsilon-\frac{3}{2}}}{\pi}\sum_{k=0}^{\infty}\frac{\Theta(k)}{k!(\lambda b+ck)^{2\upsilon}}
G_{{4} ,{2}}^{{1},{4}} \left(\ \frac{64 y^{2}}{(\lambda b+ck)^{4}}\bigg{|} \begin{array}{lll}\frac{4-2\upsilon}{4}, \frac{3-2\upsilon}{4}, \frac{2-2\upsilon}{4}, \frac{1-2\upsilon}{4};-\\0~~;~~ \frac{1}{2}\end{array} \right),
\end{eqnarray}
\begin{eqnarray}\label{RG3}
=\frac{\sqrt{y}~2^{4\upsilon}}{\pi}\sum_{k=0}^{\infty}\frac{\Theta(k)}{k!(\lambda b+ck)^{2\upsilon+1}}
G_{{4} ,{2}}^{{1},{4}} \left(\ \frac{64 y^{2}}{(\lambda b+ck)^{4}}\bigg{|} \begin{array}{lll}\frac{3-2\upsilon}{4}, \frac{2-2\upsilon}{4}, \frac{1-2\upsilon}{4}, \frac{-2\upsilon}{4};-\\\frac{-1}{4}~~;~~ \frac{1}{4} \end{array} \right),~~~~~
\end{eqnarray}
\begin{eqnarray}\label{RG4}
=\frac{2^{4\upsilon-\frac{9}{2}}}{(\pi y)}\sum_{k=0}^{\infty}\frac{\Theta(k)}{k!(\lambda b+ck)^{2\upsilon-2}}
G_{{4} ,{2}}^{{1},{4}} \left(\ \frac{64 y^{2}}{(\lambda b+ck)^{4}}\bigg{|} \begin{array}{lll}\frac{6-2\upsilon}{4}, \frac{5-2\upsilon}{4}, \frac{4-2\upsilon}{4}, \frac{3-2\upsilon}{4};-\\\frac{1}{2}~~;~~ 1 \end{array} \right),
\end{eqnarray}
\begin{eqnarray}\label{RG5}
=\frac{y~2^{4\upsilon+\frac{3}{2}}}{\pi }\sum_{k=0}^{\infty}\frac{\Theta(k)}{k!(\lambda b+ck)^{2\upsilon+2}}
G_{{4} ,{2}}^{{1},{4}} \left(\ \frac{64 y^{2}}{(\lambda b+ck)^{4}}\bigg{|} \begin{array}{lll}\frac{2-2\upsilon}{4}, \frac{1-2\upsilon}{4}, \frac{-2\upsilon}{4}, \frac{-1-2\upsilon}{4};-\\\frac{-1}{2}~~;~~0 \end{array} \right),
\end{eqnarray}
\begin{multline}\label{RG6}
=\frac{2^{6\upsilon-\frac{7}{2}}}{\pi^{3}}\sum_{k=0}^{\infty}\frac{\Theta(k)}{k!(\lambda b+ck)^{2\upsilon}}\times\\\times
G_{{8} ,{4}}^{{2},{8}} \left(\ \frac{4^{8}y^{4}}{(\lambda b+ck)^{8}}\bigg{|} \begin{array}{lll}\frac{8-2\upsilon}{8}, \frac{7-2\upsilon}{8}, \frac{6-2\upsilon}{8}, \frac{5-2\upsilon}{8}, \frac{4-2\upsilon}{8}, \frac{3-2\upsilon}{8}, \frac{2-2\upsilon}{8}, \frac{1-2\upsilon}{8};-\\ 0, \frac{1}{2}~~;~~\frac{1}{4}, \frac{3}{4} \end{array} \right).
\end{multline}
\\
\textbf{Proof.} We present the details of the proof of the right-hand side of (\ref{RG2}). Applying the Mellin-Barnes  contour integral for the cosine function in (\ref{RCG43}) and setting $x^{2}=t$, we obtain after some simplification
\begin{multline}\label{RG7}
I_{C}^{*}(\upsilon,b,c,\lambda,y)\\
=2\sqrt{\pi}\sum_{k=0}^{\infty}\frac{\Theta(k)}{k!}
\frac{1}{2\pi i}\int_{-i\infty}^{+i\infty}\frac{\Gamma(-\zeta)}{\Gamma(\frac{1}{2}+\zeta)} \left(\ \frac{y^{2}}{4}\right)^{\zeta}
\left\{\int_{0}^{\infty}e^{-(\lambda b+ck)t}~t^{2\upsilon+4\zeta-1}dt\right\}d\zeta.
\end{multline}
Evaluation of the inner integral in (\ref{RG7}) as a gamma function then leads to
 \begin{multline}\label{RG8}
I_{C}^{*}(\upsilon,b,c,\lambda,y)\\
=2\sqrt{\pi}\sum_{k=0}^{\infty}\frac{\Theta(k)}{k!(\lambda b+ck)^{2\upsilon}}
\frac{1}{2\pi i}\int_{-i\infty}^{+i\infty}\frac{\Gamma(-\zeta)\Gamma(2\upsilon+4\zeta)}{\Gamma(\frac{1}{2}+\zeta)} \left\{ \frac{y^{2}}{4(\lambda b+ck)^{4}}\right\}^{\zeta}d\zeta.
 \end{multline}
We now use the Gauss multiplication theorem (\ref{RCG50}), various algebraic  properties of Pochhammer's symbol  and the definition of the Meijer $G$-function (\ref{RCG32}) in (\ref{RG8}). After simplification we finally obtain the result in terms of an infinite series of Meijer $G$-functions stated in (\ref{RG2}).

Similarly,  the proof of the other results (\ref{RG3})--(\ref{RG6}) can be obtained
by using the Mellin-Barnes contour integrals (\ref{RCG44})--(\ref{RCG47}) in the series expansion of $I_{C}^{*}(\upsilon,b,c,\lambda,y)$ and following the same procedure as outlined above.\\
\textbf{Corollary 1}. The following infinite Fourier cosine transform of $x^{\upsilon-1}e^{-b\lambda\sqrt{x}}{}_{r}\Psi_{s}[\cdot]$ given by
\begin{equation}\label{RG9}
\Xi_{C}(\upsilon,b,c,\lambda,y)=\int_{0}^{\infty}x^{\upsilon-1}e^{-b\lambda\sqrt{x}}{}_{r}\Psi_{s}\left[\ \begin{array}{lll}(\alpha_{1}, A_{1}),...,(\alpha_{r}, A_{r});\\(\beta_{1}, B_{1}),...,(\beta_{s}, B_{s});~\end{array} e^{-c\sqrt{x}}\right]\cos (xy)\,dx,
\end{equation}
holds true for
$y>0,~\Re(\upsilon)>0,~\Re(c)>0$ and $\Re(\lambda)>0,~\Re(b)>0$ (or $\Re(\lambda)<0,~\Re(b)<0)$,
where the parameters $\alpha_{i},\beta_{j}\in\mathbb{C};$ and the coefficients
$ A_{i},B_{j}\in\mathbb{R}=(-\infty,+\infty)$, with $A_{i}\neq0~(i=1,2,...,r)$, $B_{j}\neq0~(j=1,2,...,s)$. The function ${}_{r}\Psi_{s}[\cdot]$ is the Fox-Wright function of one variable subject to suitable convergence conditions derived from those stated for (\ref{RCG21})--(\ref{RCG25}).
\\
\textbf{Proof.} The above expression is a generalized form of the Ramanujan integral in terms of the Fox-Wright function obtained by setting  \[\Theta(k)=\displaystyle{\frac{\Gamma(\alpha_{1}+kA_{1})...\Gamma(\alpha_{r}+kA_{r})} {\Gamma(\beta_{1}+kB_{1})...\Gamma(\beta_{s}+kB_{s})}},\]
for $(k=0,1,2,3,...)$ in the integral $I_{C}^{*}(\upsilon,b,c,\lambda,y)$. Then after simplification we get the result stated in (\ref{RG9}).\\
\textbf{Corollary 2.} The infinite Fourier cosine transform of $x^{\upsilon-1}e^{-b\lambda\sqrt{x}}{}_{r}F_{s}(\cdot)$ given by
\begin{equation}\label{RG10}
\nabla_{C}(\upsilon,b,c,\lambda,y)=\int_{0}^{\infty}x^{\upsilon-1}e^{-b\lambda\sqrt{x}}{}_{r}F_{s}\left(\ \begin{array}{lll}\alpha_{1},...,\alpha_{r};\\\beta_{1},...,\beta_{s};~\end{array} e^{-c\sqrt{x}}\right)\cos (xy)\,dx,
\end{equation}
holds true for $y>0,~\Re(\upsilon)>0,\Re(c)>0$ and $\Re(\lambda)>0,~\Re(b)>0$ (or $\Re(\lambda)<0,~\Re(b)<0)$,
where $\alpha_{i}\in\mathbb{C}(i=1,2,...,r),~\beta_{j}\in\mathbb{C}\setminus \mathbb{Z}_{0}^{-}(j=1,2,...,s)$ and $r\leq s+1$.\\
\textbf{Proof.} If we put the coefficients $A_{1}=...=A_{r}=1$ and $B_{1}=...=B_{s}=1$ in (\ref{RG9}), then we find after simplification the generalized expression of a Ramanujan integral denoted by $I_{C}(\upsilon,b,c,\lambda,y)$ involving the hypergeometric function ${}_{r}F_{s}(.)$ stated in (\ref{RG10}).

\section{Some applications of the generalized Ramanujan integral $\Xi_{C}(\upsilon,b,c,\lambda,y)$}
\label{Sec:3}
When $r=0,~s=1$ in the integral (\ref{RG9}), we get
\begin{equation}\label{MLF}
\int_{0}^{\infty}x^{\upsilon-1}e^{-b\lambda\sqrt{x}}~\Phi\left(B_{1},\beta_{1},e^{-c\sqrt{x}}\right)\cos (ax)\,dx.\\
\end{equation}
Here $\Phi(B_{1},\beta_{1},e^{-c\sqrt{x}})$ is known as the Wright function defined by
\begin{equation}
\Phi\left(B_{1},\beta_{1},e^{-c\sqrt{x}}\right)={}_{0}\Psi_{1}\left[\ \left.\begin{array}{lll}\overline{~~~~~~~~~}\\(\beta_{1}, B_{1})\end{array} \right| e^{-c\sqrt{x}}\right]=\sum_{k=0}^{\infty}{\frac{1} {\Gamma(\beta_{1}+B_{1}k)}\frac{\left(e^{-c\sqrt{x}}\right)^{k}}{k!}},
\end{equation}
where $\beta_{1}\in\mathbb{C},~B_{1}>0$.

When $r=0,~s=1$,~$B_{1}=\mu,~\beta_{1}=\lambda+1$ and $e^{-c\sqrt{x}}$ is replaced by $-e^{-c\sqrt{x}}$ in the integral (\ref{RG9}) then we get
\begin{equation}
\int_{0}^{\infty}x^{\upsilon-1}e^{-c\sqrt{x}}~J_{\lambda}^{\mu}(-e^{-c\sqrt{x}})\cos (ax)\,dx.
\end{equation}
Here $J_{\lambda}^{\mu}(-e^{-c\sqrt{x}})$ is known as the Wright generalized Bessel function defined by
\begin{eqnarray}
J_{\lambda}^{\mu}(-e^{-c\sqrt{x}})=\Phi\left(\mu,\lambda+1,-e^{-c\sqrt{x}}\right)={}_{0}\Psi_{1}\left[\ \left.\begin{array}{lll}\overline{~~~~~~~~~}\\(\lambda+1, \mu)\end{array}\right| -e^{-c\sqrt{x}}\right]\\
~~~~~=\sum_{k=0}^{\infty}{\frac{1} {\Gamma(\mu k+\lambda+1)}\frac{\left(-e^{-c\sqrt{x}}\right)^{k}}{k!}},
\end{eqnarray}
where $\lambda\in\mathbb{C},~\mu>0$.

When $r=1,~s=1,~\alpha_{1}=A_{1}=1$ in the integral (\ref{RG9}), we get
\begin{equation}
\int_{0}^{\infty}x^{\upsilon-1}e^{-c\sqrt{x}}~E_{B_{1},\beta_{1}}\left(e^{-c\sqrt{x}}\right) \cos (ax)\,dx.
\end{equation}
Here $E_{B_{1},\beta_{1}}(e^{-c\sqrt{x}})$ is known as the Mittag-Leffler function, defined by
\begin{eqnarray}
E_{B_{1},\beta_{1}}(e^{-c\sqrt{x}})={}_{1}\Psi_{1}\left[\ \left.\begin{array}{lll}(1,1)\\(\beta_{1}, B_{1})\end{array} \right| e^{-c\sqrt{x}}\right],~~~e^{-c\sqrt{x}}\neq0,\\
=\sum_{k=0}^{\infty}{\frac{1} {\Gamma(\beta_{1}+B_{1}k)}\frac{\left(e^{-c\sqrt{x}}\right)^{k}}{k!}},~~~~~~~~~~~~~~~~~~~~~
\end{eqnarray}
where $\beta_{1}\in\mathbb{C},~B_{1}>0$.


\section{Analytical evaluation of an infinite Fourier cosine transform of $x^{\upsilon-1}\{\exp(b\sqrt{x})-1\}^{-\lambda}$}
\label{Sec:4}
The analytical evaluation of the generalized Ramanujan integral $I_{C}(\upsilon,b,\lambda,y)$ in terms of an infinite series of Meijer $G$-functions does not appear to be available in the literature. We have the following theorem:\\
{\bf Theorem 2.}\ {\it The integral $I_{C}(\upsilon,b,\lambda,y)$ is given by
\begin{equation}\label{RG11}
I_{C}(\upsilon,b,\lambda,y) =\displaystyle{\int_{0}^{\infty}\frac{x^{\upsilon-1}\cos(xy)}{\left\{\exp(b\sqrt{x})-1\right\}^{\lambda}}dx}~~~~~~~~~~~~~~~~~~~~~~~~~~~~~~~~~~~~~~~~~~~~~~~~~~~~~~~~~~~~~~~~~~
\end{equation}
\begin{eqnarray}\label{RG12}
=\frac{2^{4\upsilon-~\frac{3}{2}}}{\pi}\sum_{k=0}^{\infty}\frac{(\lambda)_{k}}{k!(\lambda b+b k)^{2\upsilon}}
G_{{4} ,{2}}^{{1},{4}} \left(\ \frac{64 y^{2}}{(\lambda b+bk)^{4}}\bigg{|} \begin{array}{lll}\frac{4-2\upsilon}{4}, \frac{3-2\upsilon}{4}, \frac{2-2\upsilon}{4}, \frac{1-2\upsilon}{4};-\\ 0~~;~~ \frac{1}{2}\end{array} \right),
\end{eqnarray}
\begin{eqnarray}\label{RG13}
=\frac{2^{4\upsilon}\sqrt{y}}{\pi}\sum_{k=0}^{\infty}\frac{(\lambda)_{k}}{k!(\lambda b+b k)^{2\upsilon+1}}
G_{{4} ,{2}}^{{1},{4}} \left(\ \frac{64y^{2}}{(\lambda b+bk)^{4}}\bigg{|} \begin{array}{lll}\frac{3-2\upsilon}{4}, \frac{2-2\upsilon}{4}, \frac{1-2\upsilon}{4}, \frac{-2\upsilon}{4};-\\ -\frac{1}{4}~~;~~ \frac{1}{4}\end{array} \right),
\end{eqnarray}
\begin{eqnarray}\label{RG14}
=\frac{2^{4\upsilon-\frac{9}{2}}}{(\pi y)}\sum_{k=0}^{\infty}\frac{(\lambda)_{k}}{k!(\lambda b+b k)^{2\upsilon-2}}
G_{{4} ,{2}}^{{1},{4}} \left(\ \frac{64y^{2}}{(\lambda b+bk)^{4}}\bigg{|} \begin{array}{lll}\frac{6-2\upsilon}{4}, \frac{5-2\upsilon}{4}, \frac{4-2\upsilon}{4}, \frac{3-2\upsilon}{4};-\\ \frac{1}{2}~~;~~1\end{array} \right),
\end{eqnarray}
\begin{eqnarray}\label{RG15}
=\frac{y~2^{4\upsilon+\frac{3}{2}}}{\pi}\sum_{k=0}^{\infty}\frac{(\lambda)_{k}}{k!(\lambda b+b k)^{2\upsilon+2}}
G_{{4} ,{2}}^{{1},{4}} \left(\ \frac{64y^{2}}{(\lambda b+bk)^{4}}\bigg{|} \begin{array}{lll}\frac{2-2\upsilon}{4}, \frac{1-2\upsilon}{4}, \frac{-2\upsilon}{4}, \frac{-1-2\upsilon}{4};-\\ -\frac{1}{2}~~;~~0\end{array} \right),
\end{eqnarray}
\begin{multline}\label{RG16}
=\frac{2^{6\upsilon-\frac{7}{2}}}{\pi^{3}}\sum_{k=0}^{\infty}\frac{(\lambda)_{k}}{k!(\lambda b+b k)^{2\upsilon}}\times\\\times
G_{{8} ,{4}}^{{2},{8}} \left(\ \frac{4^{8}y^{4}}{(\lambda b+bk)^{8}}\bigg{|} \begin{array}{lll}\frac{8-2\upsilon}{8}, \frac{7-2\upsilon}{8}, \frac{6-2\upsilon}{8}, \frac{5-2\upsilon}{8}, \frac{4-2\upsilon}{8}, \frac{3-2\upsilon}{8}, \frac{2-2\upsilon}{8}, \frac{1-2\upsilon}{8};-\\ 0, \frac{1}{2}~~;~~\frac{1}{4}, \frac{3}{4} \end{array} \right)
\end{multline}
where $\Re(\upsilon)>0,y>0;\Re(\lambda)>0,~\Re(b)>0$.}
\\
\\
\textbf{Proof.} Use of the binomial function (\ref{RCG20}) in the integral $I_{C}^{*}(\upsilon,b,c,\lambda,y)$ shows that when $\Theta(k)=(\lambda)_{k}$ and $c=b$
\begin{equation}\label{RCG}
I_{C}(\upsilon,b,\lambda,y)=\sum_{k=0}^\infty \frac{(\lambda)_k}{k!}\int_{0}^{\infty}~x^{\upsilon-1}e^{-(\lambda b+bk)\sqrt{x}} \cos (xy)\,dx.
\end{equation}
Use of (\ref{RG2})--(\ref{RG6}) then leads after simplification to the evaluations given in (\ref{RG12}) -- (\ref{RG16}).

\section{Ramanujan's integral $\textbf{R}_C(m,n)$}
\label{Sec:5}
The evaluation of Ramanujan's integral $\textbf{R}_C(m,n)$ in terms of Meijer's $G$-function is given by the following:
\begin{eqnarray}\label{RG17}
\textbf{R}_C(m,n)=\int_{0}^{\infty}\frac{x^m \cos(\pi n x)}{\exp{(2\pi\sqrt{x})-1}}\,dx~~~~~~~~~~~~~~~~~~~~~~~~~~~~~~~~~~~~~~~~~~~~~~~~~~~~~~~~~~~~~~~~~~~~~~~~~~~~~~~~~~~~~~~~~~
\end{eqnarray}
\begin{eqnarray}\label{RG18}
=\frac{2^{4m+\frac{5}{2}}}{\pi}\sum_{k=0}^{\infty}\frac{1}{(2\pi+2\pi k)^{2m+2}}
G_{{4} ,{2}}^{{1},{4}} \left(\ \frac{64 n^{2}\pi^{2}}{(2\pi+2\pi k)^{4}}\bigg{|} \begin{array}{lll}\frac{2-2m}{4}, \frac{1-2m}{4}, \frac{-2m}{4}, \frac{-1-2m}{4};-\\ 0~~;~~ \frac{1}{2}\end{array} \right),
\end{eqnarray}
\begin{eqnarray}\label{RG19}
=\frac{2^{4m+4}\sqrt{n}}{\sqrt{\pi}}\sum_{k=0}^{\infty}\frac{1}{(2\pi+2\pi k)^{2m+3}}
G_{{4} ,{2}}^{{1},{4}} \left(\ \frac {64n^{2}\pi^{2}}{(2\pi+2\pi k)^{4}}\bigg{|} \begin{array}{lll}\frac{1-2m}{4}, \frac{-1-2m}{4}, \frac{-2m}{4}, \frac{-2-2m}{4};-\\ -\frac{1}{4}~~;~~ \frac{1}{4}\end{array} \right),
\end{eqnarray}
\begin{eqnarray}\label{RG20}
=\frac{2^{4m-\frac{1}{2}}}{n\pi^{2}}\sum_{k=0}^{\infty}\frac{1}{(2\pi+2\pi k)^{2m}}
G_{{4} ,{2}}^{{1},{4}} \left(\ \frac {64n^{2}\pi^{2}}{(2\pi+2\pi k)^{4}}\bigg{|} \begin{array}{lll}\frac{4-2m}{4}, \frac{3-2m}{4}, \frac{2-2m}{4}, \frac{1-2m}{4};-\\ \frac{1}{2}~~;~~1\end{array} \right),
\end{eqnarray}
\begin{eqnarray}\label{RG21}
=n2^{4m+\frac{11}{2}}\sum_{k=0}^{\infty}\frac{1}{(2\pi+2\pi k)^{2m+4}}
G_{{4} ,{2}}^{{1},{4}} \left(\ \frac {64n^{2}\pi^{2}}{(2\pi+2\pi k)^{4}}\bigg{|} \begin{array}{lll}\frac{-2m}{4}, \frac{-1-2m}{4}, \frac{-2-2m}{4}, \frac{-3-2m}{4};-\\ -\frac{1}{2}~~;~~0\end{array} \right),
\end{eqnarray}
\begin{multline}\label{RG22}
=\frac{2^{6m+\frac{5}{2}}}{\pi^{3}}\sum_{k=0}^{\infty}\frac{1}{(2\pi+2\pi k)^{2m+2}}\times\\\times
G_{{8} ,{4}}^{{2},{8}} \left(\ \frac {4^{8}n^{4}\pi^{4}}{(2\pi+2\pi k)^{8}}\bigg{|} \begin{array}{lll}\frac{6-2m}{8}, \frac{5-2m}{8}, \frac{4-2m}{8}, \frac{3-2m}{8}, \frac{2-2m}{8}, \frac{1-2m}{8}, \frac{-2m}{8}, \frac{-1-2m}{8};-\\ 0,~~\frac{1}{2}~~;~~\frac{1}{4},~~ \frac{3}{4} \end{array} \right),
\end{multline}
where  $m$ is a non-negative integer and $n$ is a positive rational number.\\
\textbf{Proof.} The results (\ref{RG17})--(\ref{RG22}) are obtained from (\ref{RG11})--(\ref{RG16}) upon setting $\upsilon=m+1$,$b=2\pi$, $\lambda=1$ and $y=n\pi$.

\section{Closed-form evaluations of some infinite series of $G$-functions:}
\label{Sec:6}
By putting $m=1$, $n=\frac{1}{2}, 2$, $m=n=2$ and $m=0$, $n=1,~ 2, ~4, ~6, ~\frac{1}{2},~ \frac{2}{5}$ in (\ref{RG17}) and (\ref{RG18}) and comparing with (\ref{RCG2})--(\ref{RCG4}) and (\ref{RCG10})--(\ref{RCG15}), we obtain the closed-form evaluations of the following infinite series of $G$-functions:
\begin{equation}\label{RG23}
\sum_{k=0}^{\infty}\frac{1}{(2\pi+2\pi k)^{4}}
~{G_{{4} ,{2}}^{{1},{4}}} \left( \frac{16\pi^{2}}{(2\pi+2\pi k)^{4}}\bigg{|} \begin{array}{lll} 0, -\frac{1}{4}, -\frac{1}{2}, -\frac{3}{4}\\ 0;\frac{1}{2}\end{array} \right)=\frac{\sqrt{2}(13-4\pi)}{1024\pi},~~~~~~~~~~~~~~~~~~
\end{equation}
\begin{equation}\label{RG24}
\sum_{k=0}^{\infty}\frac{1}{(2\pi+2\pi k)^{4}}
~{G_{{4} ,{2}}^{{1},{4}}} \left( \frac{64\pi^{2}}{(2\pi+2\pi k)^{4}}\bigg{|} \begin{array}{lll}0, -\frac{1}{4}, -\frac{1}{2}, -\frac{3}{4}\\ 0;\frac{1}{2}\end{array} \right)=\frac{\pi\sqrt{2}}{8192}\left(\frac{1}{2}-\frac{3}{\pi}+\frac{5}{\pi^{2}} \right),
\end{equation}
\begin{equation}\label{RG25}
\sum_{k=0}^{\infty}\frac{1}{(2\pi+2\pi k)^{6}}
~{G_{{4} ,{2}}^{{1},{4}}} \left( \frac{256\pi^{2}}{(2\pi+2\pi k)^{4}}\bigg{|} \begin{array}{lll} -\frac{1}{2}, -\frac{3}{4}, -1, -\frac{5}{4}\\ 0;\frac{1}{2}\end{array} \right)=\frac{\pi\sqrt{2}}{256\times2^{11}}\left(1-\frac{5}{\pi}+\frac{5}{\pi^{2}} \right),
\end{equation}
\begin{equation}\label{RG26}
\sum_{k=0}^{\infty}\frac{1}{(2\pi+2\pi k)^{2}}
~{G_{{3} ,{1}}^{{1},{3}}} \left( \frac{64\pi^{2}}{(2\pi+2\pi k)^{4}}\bigg{|} \begin{array}{lll} \frac{1}{4},0, \frac{-1}{4};\\ 0;\end{array} \right)=\frac{\pi(2\sqrt{2}-2)}{64},~~~~~~~~~~~~~~~~~~~~~~~~
\end{equation}
\begin{equation}\label{RG27}
\sum_{k=0}^{\infty}\frac{1}{(2\pi+2\pi k)^{2}}
~{G_{{3} ,{1}}^{{1},{3}}} \left( \frac{256\pi^{2}}{(2\pi+2\pi k)^{4}}\bigg{|} \begin{array}{lll} \frac{1}{4},0, \frac{-1}{4};\\ 0;\end{array} \right)=\frac{\pi\sqrt{2}}{128},~~~~~~~~~~~~~~~~~~~~~~~~~~~~~~~~~~~~~
\end{equation}
\begin{equation}\label{RG28}
\sum_{k=0}^{\infty}\bigg[\frac{1}{(2\pi+2\pi k)^{2}}
~{G_{{3} ,{1}}^{{1},{3}}} \left( \frac{1024\pi^{2}}{(2\pi+2\pi k)^{4}}\bigg{|} \begin{array}{lll} \frac{1}{4},0, \frac{-1}{4};\\ 0;\end{array} \right)\bigg]=\frac{\pi(3\sqrt{2}-2)}{256},~~~~~~~~~~~~~~~~~~~~~~~~
\end{equation}
\begin{equation}\label{RG29}
\sum_{k=0}^{\infty}\frac{1}{(2\pi+2\pi k)^{2}}
~{G_{{3} ,{1}}^{{1},{3}}} \left( \frac{2304\pi^{2}}{(2\pi+2\pi k)^{4}}\bigg{|} \begin{array}{lll} \frac{1}{4},0, \frac{-1}{4};\\ 0;\end{array} \right)=\frac{\pi(13\sqrt{2}-4\sqrt{6})}{1152},~~~~~~~~~~~~~~~~~~~
\end{equation}
\begin{equation}\label{RG30}
\sum_{k=0}^{\infty}\frac{1}{(2\pi+2\pi k)^{2}}
~{G_{{3} ,{1}}^{{1},{3}}} \left( \frac{16\pi^{2}}{(2\pi+2\pi k)^{4}}\bigg{|} \begin{array}{lll} \frac{1}{4},0, \frac{-1}{4};\\ 0;\end{array} \right)=\frac{\sqrt{2}}{32},~~~~~~~~~~~~~~~~~~~~~~~~~~~~~~~~~~~~~~~
\end{equation}
\begin{equation}\label{RG31}
\sum_{k=0}^{\infty}\frac{1}{(2\pi+2\pi k)^{2}}
~{G_{{3} ,{1}}^{{1},{3}}} \left( \frac{256\pi^{2}}{25(2\pi+2\pi k)^{4}}\bigg{|} \begin{array}{lll} \frac{1}{4},0, \frac{-1}{4};\\ 0;\end{array} \right)=\frac{\pi(8\sqrt{2}-3\sqrt{10})}{128}.~~~~~~~~~~~~~~~~~~
\end{equation}

\section{Conclusion}
\label{Sec:7}
We have derived an analytical expression of an infinite Fourier cosine transform related to Ramanujan's integral
as an infinite sum of Meijer $G$-functions. It is hoped that other such integrals can also be evaluated in a similar way.  We conclude by remarking that various new results and applications can be obtained from our general theorem by appropriate choice of the parameters $\upsilon,\lambda,b,c,y$ and bounded sequence $\{\Theta(k)\}_{k=0}^{\infty}$ in $I_{C}^{*}(\upsilon,b,c,\lambda,y)$.\\
\\
\textbf{Acknowledgement:}
It is good fortune and a matter of pride for me to have the esteemed guidance of Prof. R.B. Paris, Division of Computing and Mathematics, Abertay University, Dundee DD1 1HG, UK. Only his influence and valued guidance enabled me to complete the work of some research papers. He is no more now.\\
\\
\textbf{Conflict of interest:} On behalf of a present author, the corresponding author states that there is no conflict of interest.\\
\\
\textbf{Funding:} The author received no financial support for the research, authorship, and/or publication of this article.

\end{document}